\documentclass[12pt]{article}
\usepackage{MnSymbol}
\usepackage{amsmath,amsfonts,
verbatim}
\usepackage[all,cmtip]{xy}

\newcommand{\C}{{\mathbb C}}
\newcommand{\N}{{\mathbb N}}
\newcommand{\Z}{{\mathbb Z}}
\newcommand{\R}{{\mathbb R}}
\newcommand{\PP}{{\mathbf P}}
\newcommand{\HH}{\mathbb{H}}

\newcommand{\U}{{\mathbb U}}

\newcommand{\X}{{\mathbb X}}
\newcommand{\Y}{\mathbb{Y}}

\newcommand{\pr}{{\rm pr}}

\newcommand{\la}{\langle}
\newcommand{\ra}{\rangle}
\newtheorem{pkt}{}[section]  
\newcommand{\bpk}{\begin{pkt}\rm }  
\newcommand{\epk}{\end{pkt}} 
\newcommand{\inv}{^{-1}}   

\newcommand{\be}{\begin{equation}}  
\newcommand{\ee}{\end{equation}}

\newcommand{\acl}{{\rm acl}}

\newcommand{\K}{\mathrm{K}}

\newcommand{\kk}{\mathrm{k}}

\newcommand{\SL}{\mathrm{SL}}

\newcommand{\E}{\mathrm{E}}

\newcommand{\tp}{\mathrm{tp}}

\newcommand{\cl}{\mathrm{cl}}
\newcommand{\RR}{\mathsf{R}}
\newcommand{\CC}{\mathfrak{U}}
\newcommand{\card}{\mathrm{card}}

\usepackage{graphicx}

\usepackage{tikz}

\title{Non-elementary categoricity and projective locally o-minimal classes}

\author{B.Zilber}

\begin{document}
\maketitle

\abstract{Given a cover $\U$ of a family of  smooth complex algebraic varieties, we associate with it a class $\CC,$ containing $\U$, of structures locally definable  in an o-minimal expansion of the reals.  We prove that the class is $\aleph_0$-homogenous over submodels and stable. It follows that $\CC$ is categorical in cardinality $\aleph_1.$ 
In the one-dimensional case we prove that a slight modification of $\CC$ is an abstract elementary class categorical in all uncountable cardinals.

\section{Introduction}

\bpk \label{data} Let $\kk\subseteq \C,$ a subfield,  $\{ \X_i: i\in I\}$ a collection of
non-singular complex algebraic varieties defined over $\kk$ and $\U(\C),$  $\{ f_i: i\in I\},$ $I:=(I,\ge )$ a lattice with the minimal element $0$ determined by unramified $\kk$-rational epimorphisms $\pr_{i',i}:\X_{i'}\to \X_i,$ for $i'\ge i,$
a connected complex manifold $\U(\C)$ and a collection of holomorphic covering maps  (local isomorphisms)     
$$f_i: \U(\C)\twoheadrightarrow \X_i(\C),\ \ f_{i}\circ\pr_{i',i}=f_{i'}.$$
\epk
\bpk
In a number of publications, abstract elementary classes $\CC$ containing structures $(\U,f_i, \X_i),$ with an abstract algebraically closed field $\K$ instead of $\C$ (pseudo-analytic structures)  have been considered, see \cite{Zspecial} for a survey.  A typical result is a formulation of a "natural" $L_{\omega_1,\omega}$-axiom system $\Sigma$ which holds for $(\U(\C),f_i, \X_i(\C))$
 and defines a class $\CC$  categorical in all uncountable cardinals. The proofs, in each case, rely on deep results in arithmetic geometry, moreover one often is able to show that the fact of categoricity of $\Sigma$ implies the arithmetic results.
 
The above raised the question of whether an uncountably categorical AEC $\CC$ containing   $(\U(\C),f_i, \X_i(\C))$  exists under  general enough assumptions on the data, leaving aside the question of axiomatisability and related arithmetic theory.

The current paper answers this question in positive at least in the case when the $\X_i$ are curves.  We construct $\CC$
as the class of structures $\U(\K)$ ($\K=\RR+i\RR$)  locally definable (in the sense of M.Edmundo and others) in models of an o-minimal expansion $\RR$ of the reals  projected to  the language of analytic relations (we call these projective locally o-minimal). The main theorem states that, for the case when the complex dimension of $\U$ is equal to 1, an extension of $\CC$ is an abstract elementary class  categorical in all uncountable cardinals. For the general case we only where able to prove  categoricity in $\aleph_1.$

\epk

\bpk
Our main technical tool is a slightly generalised theory of $\K$-analytic sets in o-minimal expansions of the reals developed by Y.Peterzil and S.Starchenko in \cite{PZ0}. We also rely on the theory of quasi-minimal excellence as developed by the author and in the important paper \cite{K5}.  

Note that our main technical results effectively prove that the structures in $\CC$ are {\em analytic Zarsiki} in a sense slightly weaker than in the paper \cite{ZAnZar}, where we proved  results similar to the current ones for an analytic Zariski class.     

\epk
\bpk 
Most of our examples, see \ref{examples} below, have become objects of interest in the theory of o-minimality due to Pila-Wilkie-Zannier method of  counting special points of Shimura varieties and more generally, see survey \cite{Sc}. Effectively, one counts points of $\U(L)\cap D\cap S$ where $D$ is an $\RR$-definable open subset of $\U,$ $S$ a certain analytic subsets  of $\U$ and $L$
   a number field   relevant to the case at hand.
  It is easy to check, at least in the well-understood cases, that  $\U(L)$ and $S$ are definable  in the class $\CC,$ and $D$ is 0-definable in  $\RR$ (so independently on the choice of a model of the o-minimal theory).
  
   At the same time one should note that in representing $\U$ as $\U(\K),$ $\K=\RR+i\RR,$ there is a remarkable degree of freedom in the choice of a model $\RR.$ 

This raises a lot of question on the interaction between the theory of AEC and o-minimality,
 the model theory - arithmetic geometry perspective of categorical  classes and the o-minimal Pila-Wilkie-Zannier  method. 

\epk
\bpk I would like to thank Martin Bays and Andres Villaveces for some useful remarks and commentaries.
\epk
\section{Preliminaries}
\bpk 

Let $\R_{\mathrm{An}}$ be an o-minimal expansion of the reals, $\C=\R+i\R$ in the language of rings and $$\mathrm{Mod}_\mathrm{An}= \{ \RR:
\RR\equiv \R_{\mathrm{An}}\}$$
the class of models of the complete o-minimal theory. 

We write $\K$ for the algebraically closed field $\K(\RR):= \RR+i\RR.$
\epk
\bpk {\bf $( \R_\mathrm{An}, \{ f_i\})$-admissible open cover of $\U.$}

In addition to the data and notation spelled out in \ref{data},
we assume that:

\begin{itemize}
\item[(i)] There is a system of connected open subsets $D_{n}(\C)\subset \U(\C),$ $n\in \N,$ definable in $\R_{\mathrm{An}}$ without parameters,
 such that
$$\mbox{ For any  }n\in \N,\ \  
D_{n}\subseteq D_{n+1},\mbox{ and }\bigcup_n D_{n}(\C) = \U(\C);$$ 
\item[(ii)] The restriction $f_{i,n}$ of  $f_i$ on $D_{n}$ is definable in  $\R_{\mathrm{An}}$ without parameters, for each $i\in I$ and $n\in \N$ and for each $i$ there is $n$ such that
$f_i(D_n)=\X_i.$ 
\item[(iii)]
For all $i\in I,$ there is a group $\Gamma_i$ of $\RR$-definable biholomorphic transformations  on $\U,$ so that fibres of $f_i$ are $\Gamma_i$-orbits, that is 
$$f_i: \U\to \X_i\cong \U/\Gamma_i.$$
Moreover, for $i>j,$  $\Gamma_i$ is a finite index subgroup of $\Gamma_j,$ that is the cover $\pr_{i,j}: \X_i\to \X_j$ is finite.

\item[(iv)]  The system of maps $f_i, i\in I$ is {\bf $\U$-complete}: there is a chain $I_0\subseteq I$  such that
$$ \bigcap_{l\in I_0} \Gamma_l=\{ 1\}.$$

\end{itemize} 
 
\epk

\bpk\label{examples} {\bf Examples of admissible $\R_\mathrm{An}.$}

In all our examples $\R_\mathrm{An}=\R_{\exp, an},$ the reals with exponentiation and restricted analytic functions. What varies is $\U,$ $\kk$ and choice of the family $\{ f_i, D_{n}: i\in I, n\in \N\}.$

\medskip
 
1. $I=\N,$ $\U=\C,$ $f_k(z)=\exp(\frac{z}{k})$ and $D_{n}=\{ z\in \C:
-2\pi n < \mathrm{Im} z<2\pi n\}.$  

2. $\U=\C,$ $I=\N$ 
$$f_k=\exp_{\tau,k}: \C \to \E_\tau\subset \PP^2,\ \ z\mapsto \exp_{\tau}(\frac{z}{k}),$$ 
the covering map for the elliptic curve $\E_\tau$ ( $\exp_\tau$ is constructed from the Weierstrass $\mathfrak{P}$-function and its derivative $\mathfrak{P}'$, with period $k\Lambda_\tau= k\Z+\tau k\Z$). 

$D_{1}$ is  the interior of the square in $\C$ with vertices
$(-\frac{\mathrm{i}+1}{2}, -\frac{\mathrm{i}-1}{2}, \frac{\mathrm{i}+1}{2}, \frac{\mathrm{i}-1}{2}),$
and $ D_{n}= n\cdot D_{\tau,1}.$

3. $\U=\HH,$   $$D_n = \{ z \in \HH : -n/2 \le \mathrm{Re}(z) < n/2\  \& \mathrm{Im}(z) > 1/(n+1)\}.$$

For $n=1$ this is the interior of the fundamental domain of the $j$-function
$$F = \{ z \in \HH : -1/2 \le \mathrm{Re}(z) < 1/2\  \&\ \mathrm{Im}(z) > \frac{1}{2}\}$$
and the results of \cite{PZdomain} state that the restriction of $j$ to $F$ is defined in $\R_{\exp, \mathrm{an}}.$
Note that, for each $n,$ $D_n$ can be covered by finitely many shifts of $D_1$ by Moebius transformations from $\Gamma:=\SL_2(\Z).$  This allows one to define $j$ on $D_n$ in $\R_{\exp, \mathrm{an}}.$

Moreover, we can similarly consider more general functions 
$$j_N: \HH\to \Y(N)\cong\HH/\Gamma(N)$$
onto level $N$ Shimura curves. A fundamental domain for $j_N$ is a finite union of finitely many shifts of $F$ and the analysis of   \cite{PZdomain} shows that the restriction of $j_N$ on its fundamental domain is definable in $\R_{\exp, \mathrm{an}}.$ Thus we can take the family $\{ j_N\}$ to be our $\{ f_i\}$ and $\Y(N)$ to be the $\X_i.$

4.  \cite{PZdomain} supplies us with a plethora of other examples, in particular $\U=\HH_g,$ the Siegel half-space, and $\X_i$ moduli spaces of polarised algebraic varieties.   
\epk

\section{The $\K$-analytic setting}
\bpk {\bf Abstract structures definable in $\RR.$}

Let $\X_i,$ $D_{n}$ and $f_i$ defined as above in an arbitrary model $\RR$ of $\mathrm{Mod}_\mathrm{An},$ and let $\K$ be the respective algebraically closed field.
 
More precisely, we define 
$$\U(\K)=\U(\RR):=\bigcup_n D_{n}(\RR),$$
which is an $L_{\omega_1,\omega}$ interpretation of $\U$ in $\RR$ for each $i\in I.$ Now $f_i: \U(\K)\to \X_i(\K)$ is defined to be the map such that it coincides with the definable map $f_{i,\bar{n}}: D_{{n}}(\RR)\to \X_i(\K)$ for each $n\in \N.$ Note that the latter is $\K$-holomorphic in the sense of \cite{PZ0}. We will often say $\K$-holomorphic (analytic) in an extended sense: the restriction $f_{i,{n}}$ of $f_i$ to  $D_{{n}}$ is $\K$-holomorphic.

We write $D_{\bar{n}}\subset \U^m$ meaning that $\bar{n}=\la n_1,\ldots,n_m\ra\in \N^m$ and
$$D_{\bar{n}}=D_{{n_1}}\times \ldots \times D_{{n_m}}.$$
Define $f_i$ on $D_{\bar{n}}$ as $\la u_1,\ldots,u_m\ra\mapsto \la f_i(u_1),\ldots,f_i(u_m)\ra.$ This obviously extends to the map $f_i$ with the domain $\U^m.$

Note that the fundamental domains for the $f_i$ on $\U^m$ can be constructed as direct products $F_{i,k_1}\times \ldots\times F_{i,k_m},$
where $F_{i,k_j}$ are fundamental domains on $\U.$

\medskip

Below we use the notion of $\K$-analytic. 

We use notation $\kk$ for a subfield $\kk\subset \K$ such that $\kk$ contains all points definable in $\RR$ without parameters. Note that $\kk$ contains any point of the form $f_i(a)$ for $i\in I$ and a definable point $a\in D_n.$

\epk
\bpk \label{defS}
{\bf Definition.} Given $S\subset \U^m$ we say that $S$  {\bf is $L_\mathrm{AEC}(\kk)$-primitive}
 if   there are Zariski closed $Z_i\subseteq \X_i^m$  defined over $\kk,$ $i\in I_S\subseteq I,$
such that $$S= \bigcap_{i\in I_S} f_i\inv(Z_i).$$

\medskip

{\bf Remark.} Note that since $I$ is a lattice we may assume that $I_S$ is linearly ordered in $I.$

\medskip

{\bf Remark.} The equality relation is $L_\mathrm{AEC}(\kk)$-primitive.

\epk
\bpk\label{sing} {\bf Lemma.} {\em Given an $L_\mathrm{AEC}(\kk)$-primitive
$S,$ the subset $S^\mathrm{sing}$ of singular points of $S$ is also an $L_\mathrm{AEC}(\kk)$-primitive.}

{\bf Proof.} Since the $f_i$ are local isomorphisms, 
$f_i(S^\mathrm{sing})=Z_i^\mathrm{sing}.$ $\Box$
\epk
\bpk
{\bf Remark.} Note that when $\K=\C,$  $S$ is an analytic subset in the complex manifold $\U^m.$ In particular, $S\cap D_{\bar{n}}$ splits into finitely many irreducible analytic subsets since the restriction of $f_j$ to $D_{\bar{n}}$ is definable in the o-minimal structure.
 
In the general case   $S\cap D_{\bar{n}}$ is   $\K$-analytic in  $D_{\bar{n}}$ in the sense of \cite{PZ0}
  for any $i\in I$ and $\bar{n}\in \N^m.$ This is immediate by the properties of the open cover of $\U.$

\epk

 
\bpk\label{1.4}\label{1.5} {\bf Proposition.} {\em Let $S\subseteq \U^m$ be $L_\mathrm{AEC}(\kk)$-primitive  and let, fo some $n,$ $S_{j,\bar
{n}}\subseteq S\cap D_{\bar{n}}$ be a $\K$-analytic irreducible component of $S\cap D_{\bar{n}}.$ Then, for any $\bar{n}'\supseteq \bar{n}$ there is unique $S_{j,\bar
{n}'}\supseteq S_{j,\bar
{n}}$
a $\K$-analytic irreducible component of $S\cap D_{\bar{n}'}.$

The set 
$$S_j:= \bigcup_{\bar{n}'\ge \bar{n}} S_{j,\bar
{n}'}$$
 is an $L_\mathrm{AEC}(\kk')$-primitive for some algebraic extension $\kk'$ of $\kk.$

The number of such $\K$-analytic components $S_j$ of $S$ is at most countable.

 $f_i(S_j)$ is a Zariski closed $\kk'$-definable subset of $\X_i^m.$ }

{\bf Proof.} By the $\K$-analytic theory of \cite{PZ0} (see 12.1 in particular) there is a unique $\K$-analytic continuation of $S_{j,\bar
{n}}$ into  $D_{\bar{n}'},$ which gives us the first statement of Proposition.

As a consequence we get $S_j$ uniquely determined. Call it an {\bf irreducible component of $S.$}

The number of such irreducible components is at most countable since the number of components in each $D_{\bar{n}'}$ is finite.



Note that since $I$ is a lattice we may assume that $$S=\bigcap_{i\in I_0} f_i\inv(Z_i)$$ for some chain $I_0\subseteq I,$ some Zariski closed $Z_i\subseteq \X_i^m$ such that $\dim Z_i=\dim S$ and $\pr_{i,l}( Z_i)= Z_l$ for $i>l$ in $I_0.$  
 
By our argument above $S^i:=f_i\inv(Z_i)=\bigcup_{k\in K_i} S^i_k$ is a union of irreducible analytic components  with maximum dimension equal to $\dim S.$ It follows that the components of $S^i$ are also components of $S^l,$ for $i>l$ and thus  $S_j$ is a component of $f_0\inv(Z_0).$ 

We  assume without loss of generality that $Z_0$ is geometrically irreducible and
 $$S=f_0\inv(Z_0).$$ We omit the subscript $0$ below.

\medskip

Claim. $f(S_j)=Z$ and for any other component  $S_k$ of $S$ there is $\gamma\in \Gamma$ such that $\gamma\cdot S_j=S_k.$ 

 Proof.
We may assume by \ref{data} that 
 $f(D_{\bar{n}})=\X^m.$
 
By our assumption then
$$Z=f(\bigcup_{k\in K} S_k)= \bigcup_{k\in K} f(S_k\cap D_{\bar{n}})$$
where $K$ lists all the components of $S.$ But $f(S_k\cap D_{\bar{n}})$ are definable subsets of the definable $Z,$ hence by the logic compactness theorem for a big enough  finite subset $K_0\subseteq K$ 
$$Z=\bigcup_{k\in K_0} f(S_k\cap D_{\bar{n}})=\bigcup_{k\in K_0} f(S_k).$$

Let $Z^\mathrm{sing}$ the singular points of $Z$ and  $S^\mathrm{sing}$
 the singular points of $S,$ which by the fact that $f$ is a local isomorphism  are related as
 $$f\inv (Z^\mathrm{sing})=S^\mathrm{sing}.$$
By the same argument
$$\dim S_k=\dim Z=\dim S$$
for any component $S_k.$

Note that if $s\in S_j\cap S_k,$ a common point of two distinct components of $S$ then $s\in S^\mathrm{sing}.$ That is 
 $S\setminus S^\mathrm{sing},$ the analytic subset of the open set
 $\U^m\setminus S^\mathrm{sing},$ splits into non-intersecting analytic components $S_k\setminus S^\mathrm{sing}$ and we get from above
 
\be \label{ZS} Z\setminus Z^\mathrm{sing}=\bigcup_{k\in K_0} f(S_k\setminus S^\mathrm{sing}).\ee
The union on the right can not be disjoint, in fact has to be equal to one of the summands.
Indeed, consider the  fact (\ref{ZS}) in the setting of the standard model $\R_\mathrm{An},$ that is when $\K=\C$ and $\U$ is a complex manifold. Note that 
$f:\U^m\to \X^m$ is a closed map (since it is a covering map). Hence $f(S_k\setminus S^\mathrm{sing})$ are closed subsets of the connected set $Z\setminus Z^\mathrm{sing}.$ Thus $f(S_{k_0}\setminus S^\mathrm{sing})=Z\setminus Z^\mathrm{sing}$ for a $k_0\in K_0.$ 
But then, for any $k\in K_0,$ any point $s\in S_k \setminus S^\mathrm{sing}$ there is $s_0\in 
S_{k_0}\setminus S^\mathrm{sing}$ such that $f(s)=f(s_0),$ and hence 
$s=\gamma\cdot s_0$ for some $\gamma\in \Gamma.$ It follows that $S_k$ and $\gamma\cdot S_{k_0}$ intersect in a non-singular point of $S$ and thus $S_k=\gamma\cdot S_{k_0}.$
Since $K_0$ is a large enough finite subset of $K,$ we conclude that 
$$S=\bigcup_{\gamma\in \Gamma} \gamma\cdot S_j.$$ 
  Claim proved.  

\medskip

Now, for any $i\in I$ consider $$Z_{ij}:=f_i(S_j)$$
which we proved to be Zariski closed irreducible and 
$$f_i\inv(Z_{ij})= \bigcup_{\gamma\in \Gamma_i} \gamma\cdot S_j.$$  
Since by assumption $\bigcap_{l\in I_0}\Gamma_l$ is trivial, for some chain $I_0\subseteq I$ containing $i,$
 we have 
 $$S_j=\bigcap_{l\in I_0} f_l\inv(Z_{lj}).$$
 
$\Box$

\epk

\bpk {\bf Definition.} For an $m$-tuple $u$ in $\U$ and a subfield $\kk\subset \K$ the {\bf locus} of
$u$ over $\kk,$ written $\mathrm{loc}(u/\kk),$
is  
the $L_\mathrm{AEC}(\kk)$-primitive $S_u$ containing $u,$ of minimal possible dimension.



{\bf Remark.} Note that $S_u=\mathrm{loc}(u/\kk)$ is {\bf $\kk$-irreducible}, that is 
$S_u$ can not be represented as $S_1\cup S_2$ with $S_1,S_2$ $L_\mathrm{AEC}(\kk)$-primitives, both distinct from $S_u.$ This follows from the same property of Zarsiki loci for the images $f_i(u).$ 
\epk

\bpk \label{1.9} {\bf Corollary.} {\em If $\kk$ is algebraically closed then 
$\mathrm{loc}(u/\kk),$ the locus of $u$ over $\kk,$ is $\K$-analytically irreducible.}
\epk

The following trivial fact will allow us later to switch from the notion of a locus to that of a type of the form $\tp(v/\bar{u}).$

\bpk \label{newlocus} {\bf Lemma.} {\em Suppose $\bar{v}_1$ and $\bar{v}_2$ are $n$-tupes from $\U$
such that
$f_i(\bar{v}_1)=f_i(\bar{v}_2)\in \X_i(\kk)^n$ for all $i\in I.$ Then
$$\mathrm{loc}(u\bar{v}_1/\kk)=\mathrm{loc}(u\bar{v}_2/\kk).$$}

{\bf Proof.} Immediate from definitions. $\Box$

\epk

\section{$L_\mathrm{AEC}$-structures}

\bpk {\bf Definition.} Define $\CC(\RR)$ to be the 
structure with universe $\U(\RR)$ in the language of $L_\mathrm{AEC}(\kk)$-primitives.


Define $\CC$ to be the class of all  structures of the form $\CC(\RR)$ with $\mathrm{cdim}\,\RR\ge \aleph_0.$ 
\epk
\bpk \label{2.1-} {\bf Proposition.}
{\em $\CC(\RR)$ interprets over $\emptyset$  (in the first order way) the field $\K,$ points of the subfield $\kk$ and all the maps $f_i: \U\to \X_i(\K).$}

{\bf Proof.} First note that the equivalence relations on $\U$
$$E_i(u_1,u_2) :\equiv f_i(u_1)=f_i(u_2)$$
are  $L_\mathrm{AEC}(\kk)$-primitives. Thus the sets $\X_i(\K)$ are interpretable as $\U/E_i$ together with the maps $f_i: \U\to \U/E_i.$

Given a Zariski closed $Z_i\subset \X_i^m$ we have $Z_i^\U:=f_i\inv(Z_i),$ a definable subset of $\U^m.$ Thus $Z_i=f_i(Z_i^\U)$ are interpretable.

Now consider the structure  $\X(\K)_{\mathrm{Zar},\kk}$  on an infinite algebraic variety $\X(\K)$ over $\kk$ equipped with relations $Z\subseteq \X^m,$ all Zariski closed $Z$ over $\kk.$
Recall that if $\K$ is algebraically closed 
  the field structure $\K$ together with its $\kk$-points is interpretable in $\X(\K)_{\mathrm{Zar},\kk}$ (note that, for a $\kk$-rational map $g:\X(\K)\to \K,$ one can reconstruct $g(\X(\K))$ as $\X(\K)/E$
  for some equivalence relation over $\kk$).  
  
  It follows one can interpret $\K$ and $\kk$-points in $\CC(\RR)$ since we did it for $\X_i(\K)$ and the $Z_i.$ $\Box$   
 
\epk
\bpk {\bf Corollary.} {\em  Any $L_\mathrm{AEC}(\K)$-primitive is 
type-definable in  $\CC(\RR)$ using parameters.}
\epk

\bpk\label{cdim1}
Recall that an o-minimal structure $\RR$ is a pre-geometry, i.e. has a well-behaved dependence relation,  and one can define a notion of a (combinatorial) dimension $\mathrm{cdim}\,A$ of a subset $A\subseteq \RR$ (not to confuse with analytic dimension) as the cardinality of a maximal independent subset of $A.$

In particular, $\mathrm{cdim}\,\RR_0=0.$

This has the following relationship with $\dim_\RR S$ for $\K$-analytic subset $S\subseteq \RR^m$ defined over a set $C:$   


{\em Suppose $\mathrm{cdim}\,\RR/C\ge m.$ Then 
$$ \dim_\RR S\ge d \mbox{ iff there is }\la s_1,\ldots,s_m\ra\in S:\  \mathrm{cdim}(\{  s_1,\ldots,s_m\}/C) =d.$$
}

\epk
\bpk \label{types} {\bf Lemma.} {\em Suppose $$\mathrm{cdim}(\RR/\kk)\ge \aleph_0.$$
Let $\U=\U(\RR),$ 
and
$S\subset \U^m$ be an  $L_\mathrm{AEC}(\kk)$-primitive, 
$\dim_{\RR} S=d.$ Then, for any countable   family $L_{j\in J}$  of   $L_\mathrm{AEC}(\kk)$-primitives such that $\dim_{\RR} L_j<d,$ all $j\in J,$  \be\label{sl} S\setminus \bigcup_{j\in J} L_j\neq \emptyset.\ee}

{\bf Proof.} Immediate from the statement of \ref{cdim1}.$\Box$

\epk
\bpk\label{pr} {\bf Proposition} (The projection of an analytic set) {\em Let $\kk$ be algebraically closed, $T\subseteq \U^{m+1}$ be an  $L_\mathrm{AEC}(\kk)$-primitive,   and let $\pr: \U^{m+1}\to \U^m$ be the projection onto the first $m$ coordinates. Then there are an  $L_\mathrm{AEC}(\kk)$-primitive $S\subseteq \U^m,$ 
 an $i_0\in I$ and a Zariski closed  subset $R\subseteq \X_{i_0}^m$ 
defined over $\kk$ such that 
 $\dim R<\dim S$ and 
 \be \label{R}S\setminus f\inv_{i_0}(R)\subseteq \pr(T)\subseteq S\ee

 Moreover, for any $d\le \dim T-\dim S,$
 there is a Zariski closed  $R_d\subset \X_{i_0}^m$ defined over $\kk$ such that $R\subseteq R_d,$ $\dim R_d<\dim S$
and
\be \label{Rd}\pr(T)\setminus f\inv_{i_0}(R_d)=\pr_d(T)\ee
where $$\pr_d(T):= \{ s\in \pr(T): \dim  (\pr\inv(s)\cap T)\le d\}.$$
 
 }  
 
 {\bf Proof.} By \ref{1.4} $T$ is an irreducible $\K$-analytic subset of $\U^{m+1}$ and thus $$T=\mathrm{loc}(\bar{u}v/\kk)$$
 for some $\bar{u}v\in \U^{m+1},$ ($\bar{u}\in \U^m,$ $v\in \U$). 
 
 Let $$S=\mathrm{loc}(\bar{u}/\kk).$$
  By definition
$$S=\bigcap_{i\in I} f_i\inv(Z_i),\ \ T=\bigcap_{i\in I} f_i\inv(W_i)$$ 
for some Zariski closed $Z_i\subseteq \X_i^m,$ $W_i\subseteq \X_i^{m+1}$ over $\kk$ and we apply the same notation to the projection map $\pr: \X_i^{m+1}\to \X_i^m.$ We may assume that all the $Z_i$ and $W_i$ are irreducible and of dimension equal to that of $S$ and $T$ respectively.
Moreover, $$f_{i_0}(S)=Z_{i_0}\mbox{ and }f_{i_0}(T)=W_{i_0}.$$
Below, for the convenience of notation, we drop mentioning the index $i_0:$
 $$f_{}(S)=Z_{}\subseteq \X^m\mbox{ and }f_{}(T)=W_{}\subseteq \X^{m+1}.$$

Note that by the analytic continuation argument we have
$$\pr(T)\subseteq S.$$

By the basic algebraic geometry,  $\pr(W)$ is a constructible Zarsiki dense subset of $Z,$  that is there are Zariski closed $R\subset Z$ over $\kk$ such that
$Z=\pr(W)\cup R$ and $\dim R< \dim Z.$

By the basic assumptions, given arbitrary  $t\in T,$ $s=\pr( t),$
 for some $\RR$-definable open neighbourhood ${U}\subset \U^{m}$ of $s$ and open neighborhood $U\times V\subset \U^{m+1}$ of $t,$ with $V\subset \U,$
  the restriction $f_U: {U}\to \X^m$ and $f_{{U}\times V}: {U}\times V\to \X^{m+1}$ are injective.

Thus we obtain the commuting diagram with injective horizontal arrows
\be\label{dgm} \begin{array}{llll}
T\cap (U\times V)\to W\\
\ \ \ \downarrow \pr \ \ \ \ \ \ \  \ \pr \downarrow  \\
S\cap U\ \ \ \ \ \ \to \pr(W)\supseteq Z\setminus R.
\end{array}
\ee
  By comparing images of down-arrows we conclude
$$S\cap U\supseteq \pr (T\cap (U\times V))\supseteq f\inv_U (Z\setminus R) $$
 Note that $$f\inv_U (Z\setminus R)=S\cap U \setminus f\inv( R) $$
 and the choice of $R$ is independent on the choice of $U.$
 Hence $\pr(T)\supseteq S \setminus f\inv(R)$ and (\ref{R}) is proved.
 
 To prove the  second statement recall another basic fact of algebraic geometry: there is a Zariski closed $R_d\subset \X^m$ such that
 $$\pr(W)\setminus R_d=\pr_d(W):= \{ z\in \pr(W): \dim \pr\inv(z)\cap W\le d\}. $$ 
Now repeat the argument with the diagram (\ref{dgm}) with $\pr_d(W)$ in place of $\pr(W)$. This proves (\ref{Rd}). $\Box$

\epk
Recall the notion of an {\bf analytic Zariski structure}, see \cite{ZAnZar} or \cite{Zbook}.

\bpk {\bf Corollary.} {\em Under assumptions that $\kk$ is algebraically closed and $\mathrm{cdim}(\RR/\kk)\ge \aleph_0,$ the structure $\CC(\RR)$ is an  analytic Zariski structure.} 
 
 {\bf Proof.} 
 The statement of Proposition \ref{pr} asserts that the structure on $\U$ determined by $L_\mathrm{AEC}(\kk)$-primitives satisfies the key axioms (WP) and (FC) of the definition of an analytic Zariski structure. The rest of the axioms follow easily from definitions and basic algebraic geometry. $\Box$

\epk

The next statements and its proofs are similar to one of the main statements of \cite{ZAnZar} for analytic Zariski structures. More early work of M.Gavrilovich also proves this for complex analytic Zariski structures.

\bpk\label{2.2} {\bf Proposition.} {\em  $\CC$ is $\aleph_0$-homogeneous over subfields:

Suppose $\CC(\RR_1), \CC(\RR_2)\in \CC,$  $\RR_0,\RR_1,\RR_2\in \mathrm{Mod}_\mathrm{An},$ $\mathrm{cdim}(\RR_1/\RR_0)\ge \aleph_0,$
$\mathrm{cdim}(\RR_2/\RR_0)\ge \aleph_0$
and the embeddings $\RR_0\subseteq \RR_1,$ $\RR_0\subseteq \RR_1$ induce the embeddings $\CC(\RR_0)\subseteq \CC(\RR_1)$ and $\CC(\RR_0)\subseteq \CC(\RR_2).$

Let  $\kk\subseteq \K_0=\K(\RR_0)$ be algebraically closed.

Then for any
 $\bar{u}_1\in \U^m(\RR_1),$ $\bar{u}_2\in \U^m(\RR_2),$ and $w_1\in \U(\RR_1)$ such that 
$$\mathrm{loc}(\bar{u}_1/\kk)=\mathrm{loc}(\bar{u}_2/\kk)$$
there is $w_2\in \U(\RR_2)$ such that 
$$\mathrm{loc}(\bar{u}_1w_1/\kk)=\mathrm{loc}(\bar{u}_2w_2/\kk).$$
 }

{\bf Proof.} Let $S=\mathrm{loc}(\bar{u}_1/\kk)$ and $T=\mathrm{loc}(\bar{u}_1w_1/\kk).$
Note that $\bar{u}_1$ and $\bar{u}_2$ are non-singular points of $S$ and $\bar{u}_1w_1$ a non-singular point of $T,$ by \ref{sing}.


Let $d:=\dim \pr\inv(\bar{u}_1)\cap T,$ be the dimension of the fibre over $\bar{u}_1,$ and the  subset $\pr_d(T)$ as defined in \ref{pr}.
Note that by the dimension theorem of algebraic geometry $\dim \pr_d(T)=\dim S$ since $\dim \pr_d(W)=\dim S$ (in the notation of \ref{pr}).
Note also that 
$$\dim T=\dim S+d$$
since respective equality holds for the dimensions of $W$ and $Z.$ 


It follows that $\pr_d(T)$ contains all generic over $\kk'$ points of $S,$  $\bar{u}_2\in \pr_d(T)$ and thus
$$\dim \pr\inv(\bar{u}_2)\cap T=d.$$

Thus there exists $w_2$ such that $\bar{u}_2w_2\in  \pr\inv(\bar{u}_2)\cap T$ and $\dim (w_2/\bar{u}_2\kk')=d.$
Since $T$ is $\kk'$-irreducible, $$T=\mathrm{loc}(\bar{u}_2w_2/\kk').$$    
$\Box$

\epk
\bpk {\bf Lemma.} {\em Let $S\subseteq \U^{m+n}$ be an $L_\mathrm{AEC}(\kk)$-primitive and $\bar{u}\in \U^m.$ 
Let $$S_{\bar{u}}=\{ \bar{v}\in \U^n: \bar{u}\bar{v}\in S\}.$$
Then  $S_{\bar{u}}$ is an $L_\mathrm{AEC}(\kk')$-primitive, for
$\kk',$ extension of $\kk$ by co-ordinates of $f_i(\bar{u}),$ $i\in I.$ }

{\bf Proof.} By definition $S=\bigcap_{i\in I} f_i\inv(Z_i)$ for $Z_i\subseteq \X_i^{m+n}.$ 

Let, for $z_i\in \X_i^m(\K),$
$$Z_{i,z_i}=\{ x_i\in \X_i^n(\K): z_ix_i\in Z_i\}.$$  

Thus
$$S_{\bar{u}}=\{ \bar{v}\in \U^n: \bigwedge_{i\in I}f_i(\bar{u})f_i(\bar{v})\in Z_i\}=$$
$$= \bigcap_{i\in I} f_i\inv (Z_{i, f_i(\bar{u})}).$$  $\Box$
\epk
\bpk \label{overM} {\bf Corollary.} {\em $\CC$ is $\aleph_0$-homogenous over submodels:
In notations of \ref{2.2}, let $V=\U(\RR_0).$ Then, for any
$\bar{u}_1\in \U^m(\RR_1),$ $\bar{u}_2\in \U^m(\RR_1),$ $w_1\in \U^m(\RR_1)$ such that
$$\tp(\bar{u}_1/V)=\tp(\bar{u}_2/V)$$
there is  $w_2\in \U^m(\RR_2)$ such that
$$\tp(\bar{u}_1w_1/V)=\tp(\bar{u}_2w_2/V),$$
where $\tp$ is the quantifier-free type of the form (\ref{sl}). }

{\bf Proof.} Use the statement of \ref{2.2} with $\kk'=\K_0.$
\epk
\bpk {\bf Corollary.} {\em  Galois types over submodels and over $\emptyset$ in $\CC$ are quantifier free. $\CC$ is stable.}

{\bf Proof.} By \ref{types} and \ref{overM}  types have form (\ref{sl}).   
\epk

\bpk {\bf Lemma.} {\em The structure $\CC(\RR_0)$, for $\RR_0$ the prime model of the o-minimal theory,
is a prime model of $\CC.$ }


{\bf Proof.} 
An embedding $\RR_0\preccurlyeq \RR$ induces an embedding $\CC(\RR_0)\subseteq \CC(\RR).$ 
$\Box$

\epk

\bpk {\bf Theorem.} {\em Suppose $\kk$ and the language of $\R_\mathrm{An}$ is countable.

 Let $\RR_1,\RR_2\in \mathrm{Mod}_\mathrm{An}$ 
 $$\aleph_0\le \mathrm{cdim}\, \RR_1=\mathrm{cdim}\, \RR_2\le \aleph_1.$$ Then
$$\CC(\RR_1)\cong \CC(\RR_2).$$
In particular, $\CC$ is categorical in cardinality $\aleph_1.$}

{\bf Proof.}  First consider the case when  $\mathrm{cdim}\, \RR_1=\mathrm{cdim}\, \RR_2=\aleph_0.$ Then  $\CC(\RR_1)$ and  $\CC(\RR_2)$ are countable and so we can construct an isomorphism via a countable back-and-forth process  using \ref{overM}, where $\RR_0$ is the prime model of  $\mathrm{Mod}_\mathrm{An}.$ 

In case $\mathrm{cdim} \RR_1=\mathrm{cdim} \RR_2=\aleph_1$
we represent
$$\RR_1=\bigcup_{\alpha<\aleph_1} \RR_{1,\alpha}\mbox{ and } \RR_2=\bigcup_{\alpha<\aleph_1} \RR_{2,\alpha}$$
the ascending chains of elementary extensions, $\mathrm{cdim}(\RR_{i,\alpha+1}/\RR_{i,\alpha})=\aleph_0,$  for $i=1,2,$ and $\RR_{1,0}=\RR_{2,0}$ are prime models.
  Then the required isomorphism is constructed by induction on $\alpha$ using \ref{overM} with $\RR_0\cong\RR_{i,\alpha}.$ $\Box$

\epk

\section{The one-dimensional case}

\bpk
Define a closure operator $\cl: P(\U)\to P(\U)$ ($P(\U)$ the power-set of $\U$) by the condition
$$u\in \cl(\bar{w}) \mbox{ iff } \dim \mathrm{loc}(u\bar{w}/\kk)=\dim \mathrm{loc}(\bar{w}/\kk)$$ 
for $\bar{w}\subset \U$ finite. And
$$\cl(W)=\bigcup \{ \cl(\bar{w}): w\subseteq_\mathrm{fin} W\}$$
for $W$ infinite.

\epk
\bpk \label{clW} {\bf Lemma.} {\em Suppose $\cl(W)=W.$ Then the subset $f_i(W)\subset \X_i(\K)$ is closed under $\acl,$ the algebraic closure in the sense of fields. 

There is an algebraically closed subfield $L=L_W\subseteq \K.$
$$f_i(W)=\X_i(L), \mbox{ for all } i\in I.$$
} 

{\bf Proof.} Let $\bar{w}\in W^n$ and $f_i(\bar{w})=\bar{x}\in \X_i^n(\K).$ Let $y\in \X_i(\K)$ such that $y\in \acl(\bar{x}),$ where $\acl$ is over the base field $\kk.$ Thus, for
$$X=\mathrm{loc}(\bar{x}/\kk), \ Y=\mathrm{loc}(\bar{x}y/\kk)$$
we have $\dim X=\dim Y.$ Hence, since $f_i$ is a local isomorphism, for any $v\in f_i\inv(y)$ 
$$ \dim \mathrm{loc}(\bar{w}/\kk)=\dim \mathrm{loc}(\bar{w}v/\kk)$$
which implies $v\in \cl(\bar{w})\subset W.$ This proves that $f_i(W)$ is closed under $\acl$ and hence $f_i(W)=\X_i(L)$ for some algebraically closed field $L=L_{W,i}.$

We claim that $L_{W,i}=L_{W,j}$ for any $i,j\in I.$ Indeed, consider the direct product
$\U\times \U$ instead of $\U$ and $$f_i\times f_j: \U\times \U\twoheadrightarrow X_i\times X_j$$
instead of $f_i$ and $f_j,$ which still are local isomorphism onto smooth algebraic varieties. Clearly, $\cl(W\times W)=W\times W$ for $\cl$ in the product structure and
$$\X_i(L_{W,ij})\times \X_j(L_{W,ij})= (f_i\times f_j)(W\times W)=\X_i(L_{W,i})\times  \X_j(L_{W,j})$$ 
that is $L_{W,ij}=L_{W,i}=L_{W,j}=L.$ $\Box$

\epk

\bpk Recall (see \cite{K5}) that one calls $(\U,\cl)$ a  {\bf quasiminimal pregeometry structure} if the following holds:

QM1. The pregeometry is determined by the language. That is, if $\tp(v\bar{w}) = \tp(v'\bar{w}')$
 then $v\in  \cl(\bar{w})$ if and only if $v'\in  \cl(\bar{w}').$

QM2. $\U$ is infinite-dimensional with respect to $\cl.$

QM3. (Countable closure property) If $W \subset  \U$ is finite then $\cl(W)$ is countable.

QM4. (Uniqueness of the generic type) Suppose that $W, W'\subseteq \U$ are countable
subsets, $\cl(W)=W,$ $\cl(W')=W'$ and $W,W'$ enumerated so that $\tp(W) = \tp(W').$ 

If $v \in \U \setminus W$
and $v'\in  \U \setminus W'$ then $\tp(W v) = \tp(W' v')$
 (with respect to the same
enumerations for $W$ and $W'$).

QM5. ($\aleph_0$-homogeneity over closed sets and the empty set)
Let $W, W' \subseteq \U$ be countable closed subsets or empty, enumerated such
that $\tp(W) = \tp(W'),$ and let $\bar{w}, \bar{w}'$
 be finite tuples from $\U$ such that
$\tp(W\bar{w}) = \tp(W'\bar{w}')$, and let $v \in \cl(W\bar{w}).$ Then there is $v'\in \U$
that $\tp( \bar{w} vW) = \tp(\bar{w}' v'W').$

\epk
\bpk\label{QM} {\bf Proposition.} {\em Assume that $\dim_\K \U=1$ 
(the $\K$-analytic dimension) and $\mathrm{cdim}\,\RR\ge \aleph_0.$ Then $(\U(\RR),\cl)$ is a quasi-minimal pregeometry.} 

{\bf Proof.} QM1 is by definition. 

QM2 is by the assumtion on $\RR.$

QM3 follows from the well-known fact that in the language of  o-minimal structure $\RR$
$\acl(W)$ is countable and that $\cl(W)\subseteq \acl(W).$
 The latter follows from (\ref{cdim1}).
 
QM4 follows from the fact that $\U$ is one-dimensional irreducible and $v\notin \cl(W),$ $v'\notin \cl(W').$

QM5. If $W$ and $W'$ are empty then the required follows from \ref{2.2}
when $\kk'=\kk.$ In the non-empty case 
we may assume by $\aleph_0$-homogeneity over $\emptyset$ that $W=W'.$
Now  \ref{clW} and \ref{newlocus} allows to replace 
  $\tp( \bar{w} W) $ and $ \tp(\bar{w}' W')$ by 
  $\mathrm{loc}(\bar{w}/L_W)$ and  $\mathrm{loc}(\bar{w}'/L_W)$
  and $\tp( \bar{w} vW)$ and $ \tp(\bar{w}' v'W')$ by 
  $\mathrm{loc}(\bar{w}v/L_W)$ and  $\mathrm{loc}(\bar{w}'v'/L_W)$ respectively.
  
  The existence of $v'$ follows from \ref{2.2} when $\kk'=L_W.$ 
 $\Box$
\epk 
Now we recall that given a quasiminimal pregeometry structure $(\U,\cl)$
 one can associated with it an abstract elementary class containing the structure, see \cite{K5},  2.2 - 2.3, or more general \cite{ZAnZar}, 2.17 - 2.18. Call this class $\CC_\mathrm{AEC}.$ 

By definiton, one starts with a structure $\U=\CC(\RR)$ for a $\RR$ of cardinality $\aleph_1.$ Define $\CC^-_\mathrm{AEC}$ to be the class of all $\cl$-closed substructures of $\U$ with embedding $\prec$ of structures defined as {\em closed} embeding, that is

$\U_1\prec \U_2$ if and only if $\U_1\subset \U_2$ and, for finite $W\subset \U_1,$ $$\cl_{\U_1}(W)=\cl_{\U_2}(W).$$

Now define  $\CC_\mathrm{AEC}$ to be the smallest class which contains
$\CC^-_\mathrm{AEC}$ and is closed under unions of  $\prec$-chains.

\bpk\label{cc} {\bf Lemma.}
 $$\CC\subseteq  \CC_\mathrm{AEC}.$$
 
 {\bf Proof.} We need to show that $\U(\RR)\in \CC_\mathrm{AEC},$ for any $\RR\in \mathrm{Mod}_\mathrm{An}.$

We prove by induction on  $\kappa=\card\, \RR\ge \aleph_1$
 that there is a $\kappa$-chain $$\{ \U_\lambda\in \CC_\mathrm{AEC}: \lambda\in\kappa\}\mbox{ such that } \bigcup_{\lambda\in\kappa} \U_\lambda=\U(\RR).$$

Indeed,  $\RR$ can be represented as
$$\RR=\bigcup_{\lambda<\kappa} \RR_\lambda$$
for an elementary chain 
$$\{ \RR_\lambda\in : \lambda\in\kappa\}, \ \card\,\RR_\lambda=\card\, \lambda+\aleph_0, \ \RR_\lambda\prec \RR_\mu\mbox{ for } \lambda<\mu.$$
Hence
$$\U_\lambda:=\CC(\RR_\lambda)\in \CC_\mathrm{AEC}$$
which proves the inductive step and the lemma. $\Box$ 
 
\epk    
\bpk {\bf Theorem.} {\em Assuming  $\dim_\K \U=1,$ the class $\CC_\mathrm{AEC}$ is an abstract elementary class extending $\CC.$
  $\CC_\mathrm{AEC}$ is categorical in uncountable cardinals and can be axiomatised by an $L_{\omega_1,\omega}(Q)$-sentence.}

{\bf Proof.} The first part is by \ref{QM} - \ref{cc} above. 
The second part is the main result, Theorem 2.3, of \cite{K5}. $\Box$  
\epk

\end{document}